\documentclass{amsart}
\usepackage{amsmath}
\usepackage{amssymb}
\usepackage{amsfonts}
\usepackage{amsthm}
\usepackage{enumerate}
\usepackage[all]{xy}
\usepackage[latin1]{inputenc}
\usepackage{mathdots}
\usepackage{graphicx}
\usepackage{latexsym}
\usepackage{mathrsfs}
\usepackage{textcomp}
\usepackage{color}
\input xy

\makeatletter

\newtheoremstyle{definition}
        {5pt}
        {3pt}
        {}
        {0pt}
        {\scshape}
        {.}
        {5pt}
        {\thmname{#1} \thmnumber{#2} \thmnote{[#3]}} 

\newtheoremstyle{theorems}
        {5pt}
        {3pt}
        {\itshape}
        {0pt}
        {\scshape}
        {.}
        {5pt}
        {\thmname{#1} \thmnumber{#2}\thmnote{[#3]}} 

\swapnumbers

\theoremstyle{theorems}
\newtheorem{Theo}{Theorem}[section]
\newtheorem{Prop}[Theo]{Proposition}

\newtheorem{Lemma}[Theo]{Lemma}

\theoremstyle{definition}
\newtheorem{Defn}[Theo]{Definition}

\newcommand{\Oa}{\mathit\Omega}
\newcommand{\Sa}{\mathit\Sigma}
\newcommand{\Ga}{\mathit\Gamma}
\newcommand{\Da}{\mathit\Delta}
\newcommand{\cC}{{\mathcal{C}}}
\newcommand{\GaA}{\Ga\hspace{-2pt}_A}

\newcommand{\Hom}{{\rm Hom}}

\newcommand{\mmod}{{\rm mod}\hspace{0.4pt}}
\newcommand{\ind}{{\rm ind}\hspace{0.4pt}}
\newcommand{\rad}{{\rm rad}}

\newcommand{\taum}{{\tau^-\hspace{-1pt}}}

\makeatother
\newcommand{\A}{{\mathbb{A}}}
\newcommand{\Z}{{\mathbb{Z}}}
\newcommand{\N}{{\mathbb{N}}}

\begin{document}

\title{\sc Auslander-Reiten components with bounded short cycles}

\keywords{Representation-finite algebras; tilted algebras; Jacobson radical; cycles in a module category; irreducible maps; Auslander-Reiten sequences; translation quivers; Auslander-Reiten quiver; Auslander-Reiten components.}

\subjclass[2010]{16G30, 16G60, 16G70}

\author[Shiping Liu]{Shiping Liu}

\address{Shiping Liu\\ D\'epartement de math\'ematiques, Universit\'e de Sherbrooke, Sherbrooke, Qu\'ebec, Canada}
\email{shiping.liu@usherbrooke.ca}

\author[Jinde Xu]{Jinde Xu}

\address{Jinde Xu\\ D\'epartement de math\'ematiques, Universit\'e de Sherbrooke, Sherbrooke, Qu\'ebec, Canada}
\email{jinde.xu@usherbrooke.ca}

\thanks{The first author is supported in part by the Natural Sciences and Engineering Research Council of Canada.}

\maketitle

\begin{abstract}

\vspace{-20pt}

We study Auslander-Reiten components of an artin algebra with bounded short cycles, namely, there exists a bound for the depths of maps appearing on short cycles of non-zero non-invertible maps between modules in the given component. First, we give a number of combinatorial characteri\-zations of almost acyclic Auslander-Reiten components. Then, we show that an Auslander-Reiten component with bounded short cycles is closely related to the connec\-ting component of a tilted quotient algebra. In particular, the number of such components is finite and each of them is almost acyclic with only finitely many DTr-orbits. As an application, we show that an artin algebra is representation-finite if and only if its module category has bounded short cycles. This includes a well known result of Ringel's, saying that a representation-directed algebra is representation-finite.

\end{abstract}

\bigskip

\section*{Introduction}

\medskip

Let $A$ be an artin algebra, and let ${\rm mod}\hspace{.5pt}A$ stand for the category of finitely generated left $A$-modules. We shall denote by ${\rm rad}({\rm mod}\hspace{.5pt}A)$ and $\GaA$ the Jacobson radical and the Auslander-Reiten quiver of ${\rm mod}\hspace{.5pt}A$, respectively. The ultimate objective of representation theory is to study the full subcategory ${\rm ind}\hspace{.5pt}A$ of ${\rm mod}\hspace{.5pt}A$ generated by the indecomposable mo\-dules. In \cite{Rin}, Ringel initiated the study of directing modules, that is, indecomposable mo\-dules not lying on any cycle of non-zero non-invertible maps in ${\rm ind}\hspace{.5pt}A$. He showed that directing modules are uniquely determined (up to isomorphism) by their composition factors, and $A$ is representation-finite if ${\rm ind}\hspace{.5pt}A$ has no cycle; see \cite[(2.4)]{Rin}. In general, the directing mo\-dules fall into finitely many DTr-orbits of $\GaA$; see \cite{PeX, SkS}. Later, these results have been gene\-ralized to indecomposable modules not lying on any short cycle (that is, cycles of at most two maps) in ${\rm ind}\hspace{.5pt}A$. For instance, these modules are also uniquely determined (up to isomorphism) by their composition factors; see \cite{RSS1}, and $A$ is representation-finite if ${\rm ind}\hspace{.5pt}A$ has no short cycle; see \cite{HL}. Moreover,  $\GaA$ has at most finitely connected components $\Ga$ such that ${\rm rad}(\Ga)$ has no short cycle, and each of such components contains only finitely many DTr-orbits; see \cite[(2.6),(2.8)]{Liu4}.

\medskip

On the other hand, the representation theory of $A$ is determined to certain extend by the Jacobson radical of ${\rm mod}\hspace{.5pt}A$. Indeed, a well known result of Auslander's says that $A$ is representation-finite if and only if the infinite radical of ${\rm mod}\hspace{.5pt}A$ vanishes, or equivalently, every non-zero map in ${\rm mod}\hspace{.5pt}A$ is of finite depth; see \cite[(1.1)]{Sko}. For many best understood classes of representation-infinite algebras, such as tame concealed algebras, the cycles in their module category contains only maps of finite depth; see \cite{Rin}. Motivated by this fact, Skowro\'nski studied in \cite{Sko4} cycles of maps of finite depth, which are originally called {\it finite cycles}.
Indeed, Auslander-Reiten components whose non-directing modules lie only on cycles of maps of finite depth are described in \cite{MPS2}, and algebras whose module category contains only cycles of maps of finite depth are extensively studied by many authors; see, for example, \cite{MPS1, MPS2, Sky, Sko4}. More gene\-rally, the connected components $\Ga$ of $\GaA$ such that ${\rm rad}(\Ga)$ contains only short cycles of maps of finite depth are studied in \cite{Liu4}.

\medskip

The main purpose of this paper is to investigate the connected components $\Ga$ of $\GaA$ such that ${\rm rad}(\Ga)$ has bounded short cycles, that is, there exists a bound for the depths of the maps appearing on short cycles in ${\rm rad}(\Ga)$. To start with, we shall study almost acyclic Auslander-Reiten components, which have played a special role in the study of generalized double tilted algebras introduced by Reiten and Skowro\'nski; see \cite{RS}. Our main result says that if $\Ga$ is an infinite semi-stable component of $\GaA$ with bounded short cycles, then $A/{\rm ann}(\Ga)$ is a tilted algebra whose connecting component contains most of the modules in $\Ga$; see (\ref{scbc-ta}). In particular, given a connected component $\cC$ of $\GaA$, if $\cC$ is semi-regular, then ${\rm rad}(\cC)$ has bounded short cycles if and only if $A/{\rm ann}(\cC)$ is tilted with $\cC$ being its connecting component; see (\ref{semi-reg}); and if $\cC$ is generalized standard, then ${\rm rad}(\cC)$ has bounded short cycles if and only if $A/{\rm ann}(\cC)$ is generalized double tilted algebra with $\cC$ being its connecting component; see (\ref{main-th-3}). In general, $\GaA$ has at most finitely many connected components with bounded short cycles, and each of them is almost acyclic with only finitely many DTr-orbits; see (\ref{almost_acyclic}).
Finally, we shall show that an artin algebra is representation-finite  if and only if ${\rm ind}\,A$ has bounded short cycles; see (\ref{last-thm}). This includes the above-mentioned result of Ringel's as a special case.

\smallskip

\section{Preliminaries}

\medskip

\noindent The objective of this section is to collect some terminology and fix some notation which will be used throughout this paper.


\subsection{\sc Quivers.} All quivers in this paper are locally finite. Let $Q$ be a quiver. We shall say that $Q$ is {\it trivial} if it consists of a single vertex, {\it acyclic} if $Q$ contains no oriented cycle, and {\it almost acyclic} if $Q$ contains at most finitely many vertices which lie on some oriented cycles. Given two vertices $a, b$, the {\it interval} $[a, b]$ is the set of vertices lying on some path from $a$ to $b$ in $Q$. We shall say that $Q$ is {\it interval-finite} if all intervals in $Q$ are finite; and {\it strongly interval-finite} if, given any vertices $a, b$, the number of paths in $Q$ from $a$ to $b$ is finite. If $Q$ is strongly interval-finite, then it is clearly interval-finite and acyclic. A full subquiver $\Sa$ of $Q$ is called {\it convex in $Q$} if it contains all the paths with end-points lying in $\Sa$; {\it predecessor-closed} if it contains all the predecessors of its vertices in $Q$; and {\it successor-closed} if it contains all the successors of its vertices in $Q$. Finally, an infinite path in $Q$ is called {\it left infinite} if it has no starting point; {\it right infinite} if it has no ending point; and {\it double infinite} if it has neither starting point nor ending point.

\medskip

\subsection{\sc Translation Quivers.} Let $\Ga$ be a translation quiver with translation $\tau$; see, for definition, \cite[Page 47]{Rin}. A path $\hspace{-2pt}\xymatrixcolsep{18pt}\xymatrix{a_0\ar[r] & a_1\ar[r]& \cdots \ar[r] & a_n}$ in $\Ga$ is called {\it sectional} if there exists no $0<i<n$ such that $a_{i-1}=\tau a_{i+1}.$ An infinite path is {\it sectional} if all its finite subpaths are sectional. A vertex $a$ of $\Ga$ is called {\it left stable} (respectively, {\it right stable}) if $\tau^na$ is defined for all $n\ge 0$ (respectively, for all $n\le 0$); {\it stable} if its left and right stable; and $\tau$-{\it periodic} if $\tau^na=a$ for some $n>0$. Moreover, $\Ga$ is called {\it left stable} (respectively, {\it right stable}, {\it stable}, $\tau$-{\it periodic}) if every vertex of $\Ga$ is left stable (respectively, right stable, stable, $\tau$-periodic).

\medskip

The full subquivers $_l\Ga$, $_r\Ga$ and $_s\Ga$ of $\Ga$ generated by the left stable vertices, by the right stable vertices, and by the stable vertices are called the {\it left stable part}, the {\it right stable part}, and the {\it stable part} of $\Ga$ respectively. Moreover, the connected components of the quiver $_l\Ga$ (respectively, $_r\Ga$, $_s\Ga$) are called the {\it left stable components} (respectively, {\it right stable components, stable components}) of $\Ga$; and a left or right stable component of $\Ga$ is simply called a {\it semi-stable component}. Since $\Ga$ is locally finite, a semi-stable component of $\Ga$ is $\tau$-periodic if it contains a $\tau$-periodic vertex.

\medskip

Given an acyclic quiver $\Da$, one constructs in a canonical way a stable translation quiver $\Z\Da$ with translation $\rho$; see, for example, \cite[Section 2]{Liu}.
If $\Da$ is of type $\A_\infty$, then $\mathbb{Z}\Da$ does not depend on the orientation of $\Da$, 
and hence, $\Z\Da$ will be simply written as $\Z\A_\infty$. A translation quiver is called a {\it stable tube} if it is isomorphic to $\Z \A_\infty/\hspace{-3pt}<\hspace{-2pt}\rho^n\hspace{-2pt}>$ for some integer $n>0$; and {\it quasi-serial} if it is a stable tube or of shape $\Z\A_\infty$. Starting with a quasi-serial translation quiver, one obtains new translations by {\it ray insertions} or by {\it co-ray insertions}; see \cite[Section 2]{Liu1}.

\medskip

Let $\Sa$ be a connected full subquiver of $\Ga$. Recall that $\Sa$ is a {\it section} of $\Ga$ if it is acyclic, convex, and contains exactly one vertex of each $\tau$-orbit in $\Ga$; see \cite[(2.1)]{Liu}. In this case,
there exists an embedding $\Ga\to \Z \Sa$, which sends a vertex $\tau^nx$ with $n\in \Z$ and $x\in \Sa$ to the vertex $(-n, x)$; see \cite[(2.3)]{Liu}.
More generally, 
$\Sa$ is called a {\it cut} of $\Ga$; see \cite[(2.1)]{Liu5} provided, for any arrow $a\to b$ in $\Ga$, that the following two conditions are verified.

\vspace{-1.5pt}

\begin{enumerate}[$(1)$]

\item If $a\in \Sa$, then either $b$ or $\tau b$, but not both, lies in $\Sa$;

\vspace{.5pt}

\item If $b\in \Sa$, then either $a$ or $\tau^-a$, but not both, lies in $\Sa$.

\end{enumerate}


\subsection{\sc Module Category.} Throughout this paper, $A$ stands for an artin algebra over a commutative artinian ring $R$. We shall denote by $\mmod A$ the category of finitely generated left $A$-modules, and by $\rad(\mmod A)$ the Jacobson radical of $\mmod A$. Recall that $\rad^\infty(\mmod A)=\cap _{\,n\ge 0}\, \rad^n(\mmod A)$, where $\rad^n(\mmod A)$ stands for the $n$-th power of $\rad(\mmod A)$, is called the {\it infinite radical} of $\mmod A$. Given a map $f$ in $\mmod A$, its {\it depth} ${\rm dp}(f)$ is $\infty$ if $f\in {\rm rad}^\infty(X,Y)$; and otherwise, the minimal integer $n\geq 0$ for which $f\in \rad^n(X,Y) \backslash \rad^{n+1}(X,Y)$; see \cite[(1.2)]{Liu5}.

\medskip

We shall denote by ${\rm ind}\hspace{0.5pt}A$ a full subcategory of $\mmod A$ generated by a complete set of representatives of the isomorphism classes of the indecomposable modules in $\mmod A$. A {\it path} of length $n$ in ${\rm ind}\hspace{.5pt}A$ is sequence  $$\xymatrix{\sigma: \; X_0\ar[r]^-{f_1} & X_1  \ar[r] & \cdots \ar[r] & X_{n-1}\ar[r]^-{f_n} & X_n}$$ of $n$ non-zero maps in ${\rm rad}({\rm ind}\hspace{.5pt}A)$, whose {\it depth} is defined to be the supremum of the ${\rm dp}(f_i)$ with $1\le i\le n$ and written as ${\rm dp}(\sigma)$. Such a path $\sigma$ is called a {\it cycle} if $X_n=X_0$; and a {\it short cycle} if, in addition, $n\le 2$.


\subsection{\sc Auslander-Reiten Quiver.} The reader is referred to \cite{ARS} for the Auslander-Reiten theory of irreducible maps and almost split sequences in $\mmod A$. For each module $X$ in $\ind\,A$, we set $D_X={\rm End}(X)/\rad({\rm End}(X))$. The {\it Auslander-Reiten quiver} $\GaA$ of $A$ is a valued translation quiver defined as follows; see, for example, \cite[Section 2]{Liu2}. The vertices are the objects of ${\rm ind}\hspace{.5pt}A$. Given two vertices $X, Y$, there exists an arrow $X\to Y$ in $\GaA$ with valuation $(d_{XY}, d'_{XY})$ if and only if there exists an irreducible map $f: X\to Y$ in $\mmod A$, where $d_{XY}$ and $d'_{XY}$ are the dimensions of $\rad(X, Y)/\rad^2(X, Y)$ over $D_Y$ and over $D_X$, respectively. The translation $\tau$ is the Auslander-Reiten translation DTr so that $\tau Z=X$ if and only if there exists has an almost split sequence $\xymatrixcolsep{18pt}\xymatrix{0\ar[r]& X \ar[r] & Y \ar[r] & Z \ar[r] & 0}$ in $\mmod A$. The connected components of $\GaA$ are called {\it Auslander-Reiten components} of $A$, and a connected component of $\GaA$ is called {\it semi-regular} if it contains no projective module or no injective module; see \cite{Liu0}.

\medskip

Let $\cC$ be a connected component of $\GaA$ with a full subquiver $\Sa$. Given $n\in \Z$, we shall denote by $\tau^n\Sa$ the full subquiver (possibly empty) of $\cC$ generated by the modules $\tau^nX$ with $X\in \Sa$. The {\it annihilator} of $\Sa$ is ${\rm ann}(\Sa)=\cap _{X\in \Sa}\; {\rm ann}(X)$, where ${\rm ann}(X)$ stands for the annihilator of $X$ in $A$. One says that $\Sa$ is {\it faithful} if ${\rm ann}(\Sa)=0$; {\it sincere} if every simple $A$-module is a composition factor of some module in $\Sa$;  {\it convex in ${\rm ind}\hspace{.5pt}A$} if every path in ${\rm ind}\hspace{.5pt}A$ with end-points lying in $\Sa$ passes only modules in $\Sa$; and {\it generalized standard} if $\rad^{\infty}(X,Y)=0$ for all modules $X, Y\in \Sa;$ compare \cite{Sko3}. Recall that $\Sa$ is a {\it slice} in $\mmod A$ if $\Sa$ is a cut of $\mathcal{C}$ which is sincere and convex in ${\rm ind}\hspace{.5pt}A$; see \cite{Rin}, and in this case, $\mathcal{C}$ is called a {\it connecting component} of $\GaA$. Finally, $A$ is {\it tilted} if $A={\rm End}_H(T),$ where $H$ is a hereditary artin algebra and $T$ is a tilting module in ${\rm mod}\hspace{0.5pt}H;$ see \cite{HaR}. If $A$ is connected, then it tilted if and only if $\GaA$ has a connecting component; see \cite{Rin}.

\medskip

\section{Almost acyclic components}

\medskip

The objective of this section is to study some combinatorial properties of the Auslander-Reiten quiver $\GaA$. This will allow us to obtain a number of combinatorial characterizations of the almost acyclic components of $\GaA$.
Let us start with an easy observation.

\medskip

\begin{Lemma}\label{num-sec-path}

If $X, Y$ are modules in $\GaA$, then the number of sectional paths in $\GaA$ from $X$ to $Y$ is finite.

\end{Lemma}

\noindent{\it Proof.} Let $X, Y$ be modules in $\GaA$. Since $\Hom_A(X, Y)$ is of finite $R$-length, there exists some integer $r>0$ such that
$\rad^r(X, Y)=\rad^\infty(X, Y)$. Suppose that $$\xymatrix{X=X_0\ar[r] & X_1\ar[r]& \cdots \ar[r] & X_{s-1}\ar[r] &
X_s=Y}$$ is a sectional path in $\GaA$. Choosing irreducible maps $f_i: X_{i-1}\to X_i$ in $\mmod A$, $i=1, \ldots, s$, we obtain a map $f_s\cdots f_1$ of depth $s$; see \cite[(13.3)]{IgT}. In particular, $s<r$. Since $\GaA$ is locally finite, by K\"{o}nig's Lemma, the number of paths from $X$ to $Y$ of length at most $r$ is finite. The proof of the lemma is completed.

\medskip

Next, we shall study some properties of the semi-stable components of $\GaA$.

\medskip

\begin{Lemma} \label{str-int-fin}

Let $\Ga$ be a semi-stable component of $\GaA$. If $\Ga$ is acyclic, then it is strongly interval-finite.

\end{Lemma}

\noindent{\it Proof.} We shall consider only the case where $\Ga$ is a non-trivial acyclic left stable component. Observe that $\Ga$ contains a section $\Da$; see \cite[(3.3)]{Liu1}. Let $X, Y\in \Ga$ be such that $\Ga$ contains some paths from $X$ to $Y$. Since $\Ga$ embeds in $\Z\Da$, we see that $X\in \tau^s\Da$ and $Y\in \tau^t\Da$ with $s\ge t$. We shall proceed by induction on $t-s$. Assume that $s=t$. Since $\Da$ is convex in $\Ga$, so is $\tau^s\Da$. Thus, all the paths $X\rightsquigarrow Y$ in $\Ga$ belong to $\Da$, and in particular, they are all sectional. By Lemma \ref{num-sec-path}, the number of paths in $\Ga$ from $X$ to $Y$ is finite.

Suppose that $s>t$ but $\Ga$ contains infinitely many paths $\rho_n: X \rightsquigarrow Y,$ $n\in \N$. Since $\Ga$ is locally finite, we may assume that the lengths $\ell(\rho_n)$ are pairwise different. By Lemma \ref{num-sec-path}, we may assume further that none of the $\rho_n$ is sectional. Then, each path $\rho_n$ induces a path $\sigma_n: X\rightsquigarrow \tau Y$ in $\Ga$ of length $\ell(\rho_n)-2$. In particular, the $\sigma_n$ with $n\in \N$ are pairwise different paths in $\Ga$ from $X$ to $\tau Y$. Since $\tau Y\in \tau^{t+1}\Da$, we obtain a contradiction to the induction hypothesis. The proof of the lemma is completed.

%

\medskip

\begin{Lemma}\label{ss-cpt}

Let $\Ga$ be a semi-stable component of $\GaA$, not being $\tau$-periodic. If $\Ga$ has oriented cycles, then it contains a module $M$ with two infinite sectional paths
$$\xymatrixcolsep{18pt}\xymatrix{
\cdots \ar[r] & \tau^{2t}M_1\ar[r] &  \tau^t M_r\ar[r] & \cdots \ar[r] & \tau^tM_1\ar[r] &
M_r \ar[r] & \cdots \ar[r] & M_1=M
}\vspace{-1pt}$$ and
$$\xymatrixcolsep{18pt}\xymatrix{
M=N_1\ar[r]& \cdots \ar[r] & N_s\ar[r] & \tau^tN_1\ar[r] & \cdots \ar[r] & \tau^tN_s\ar[r] & \tau^{2t}N_1\ar[r] & \cdots
}$$
where $r, s, t$ are some positive integers.

\end{Lemma}

\noindent{\it Proof.} We shall consider only the case where $\Ga$ is a left stable component, which contains oriented cycles. Having no $\tau$-periodic module, $\Ga$ contains a sectional path
$$\xymatrix{
\tau^m\hspace{-1.5pt} X_1 \ar[r] & X_n\ar[r] & \cdots\ar[r] &
X_1,}$$
where $m>n\ge 1,$ and $X_1, \ldots, X_n$ lie in pairwise different $\tau$-orbits; see \cite[(2.2)]{Liu0}. This yields two sectional paths in $\Ga$ as follows:
$$\xymatrixcolsep{20pt}\xymatrix{\tau^m\hspace{-1.5pt} X_n \ar[r] & \tau^m X_{n-1}\ar[r] & \cdots \ar[r] & \tau^m X_1\ar[r] &
X_n}$$ and $$\xymatrixcolsep{20pt}\xymatrix{
X_n\ar[r] &\tau^{m-1}\hspace{-1.5pt} X_1 \ar[r] & \cdots \ar[r] & \tau^{m-(n-1)} X_{n-1} \ar[r] & \tau^{m-n} X_n.}$$

Writing $Y_1=Z_1=X_n$ and setting $Y_i=\tau^mX_{i-1}$ and $Z_i=\tau^{m-(i-1)}X_{i-1},$ for $i=2, \ldots, n$, we obtain two sectional paths $\xymatrixcolsep{18pt}\xymatrix{\tau^m Y_1 \ar[r] & Y_n \ar[r] & \cdots \ar[r] & Y_1}$ and $\xymatrixcolsep{18pt}\xymatrix{
Z_1 \ar[r] & \cdots \ar[r] & Z_n \ar[r] & \tau^{m-n}\hspace{-1.5pt}Z_1.}$
Applying repeatedly $\tau^m$ and $\tau^{m-n}$ to them, we get two infinite sectional paths
$$\xymatrixcolsep{20pt}\xymatrix{
\cdots\ar[r]&\tau^{2m} Y_1 \ar[r] & \tau^mY_n\ar[r] & \cdots \ar[r] & \tau^m Y_1 \ar[r] & Y_n\ar[r] & \cdots \ar[r] & Y_1}$$ and
$$\xymatrixcolsep{20pt}\xymatrix{
Z_1 \ar[r] & \cdots \ar[r] & Z_n\ar[r] & \tau^{m-n}\hspace{-1.5pt}Z_1\ar[r] & \cdots \ar[r] &
\tau^{m-n}\hspace{-1.5pt}Z_n\ar[r] & \tau^{2(m-n)}\hspace{-1.5pt}Z_1 \ar[r] & \cdots } $$
in $\Ga$. Writing $t=m(m-n)$ and renaming the modules in the above infinite paths so that $M_1=Y_1$ and $N_1=Z_1$, we obtain
two desired infinite sectional paths
$$\xymatrixcolsep{20pt}\xymatrix{
\cdots\ar[r]&\tau^{2t}\hspace{-1.5pt}M_1 \ar[r] & \tau^t\hspace{-1.5pt}M_r\ar[r] & \cdots \ar[r] & \tau^t\hspace{-1.5pt}M_1 \ar[r] & M_r\ar[r] & \cdots \ar[r] & M_1}$$ and
$$\xymatrixcolsep{20pt}\xymatrix{
N_1\ar[r] & \cdots \ar[r] & N_s \ar[r] & \tau^t \hspace{-1.5pt} N_1 \ar[r] & \cdots \ar[r] &
\tau^t\hspace{-1.5pt}N_s\ar[r] & \tau^{2t}\hspace{-1.5pt}N_1 \ar[r] & \cdots } $$
with $M_1=N_1$. The proof of the lemma is completed.

\medskip

%
%
%
%

We now state some properties of semi-stable but not $\tau$-periodic modules.

\medskip

\begin{Lemma}\label{ss-module}

Let $X$ be a module in $\GaA$ which is not $\tau$-periodic.

\vspace{-1.5pt}

\begin{enumerate}[$(1)$]

\item If $X$ is left stable, then there exists some $s\ge 0$ such that the predecessors of $\tau^sX$ in $\mathit\Gamma_A$ all are left stable.

\item If $X$ is right stable, then there exists some $t\ge 0$ such that the successors of $\tau^{-t}X$ in $\GaA$ all are right stable.

\end{enumerate} \end{Lemma}

\noindent{\it Proof.} We shall consider only the case where $X$ is left stable. Since $X$ is not $\tau$-periodic, there exists some $r\ge 0$ such that none of the $\tau^iX$ with $i\ge r$ has an immediate projective predecessor in $\GaA$. Then, the $\tau^iX$ with $i\ge r$ belong to a non-trivial left stable component $\Ga$ of $\GaA$. Since $X$ is not $\tau$-periodic, $\Ga$ contains no $\tau$-periodic module.

We claim that $\Ga$ contains a connected subquiver $\Sa$ such that the modules in $\Sa$ form a complete set of $\tau$-orbit representatives of $\Ga$ and have no projective predecessor in $\GaA$. Indeed, if $\Ga$ contains no oriented cycle, then it contains a section $\Sa$ with the claimed property; see \cite[(3.3)]{Liu1}. Otherwise, $\Ga$ contains a path
$$\xymatrixcolsep{20pt}\xymatrix{\tau^mX_1\ar[r] & X_n\ar[r] & \cdots \ar[r] & X_1,}$$ where $m>n\ge 1,$ and $X_1, \ldots, X_n$ form a complete set of representatives of the $\tau$-orbit in $\Ga$; see \cite[(3.6)]{Liu1}. For each $1\le j\le n$, since $X_j$ is not $\tau$-periodic, there exists some $s_j\ge 0$ such that none of the $\tau^iX_j$ with $i\ge s_j$ has a projective immediate predecessor in $\GaA$. Setting $t={\rm max}\{s_1, \ldots, s_n\}$, we see easily that
$$\xymatrixcolsep{20pt}\xymatrix{\Sa: \, \tau^tX_n\ar[r] & \cdots \ar[r] & \tau^tX_1}$$
has the claimed properties. This establishes our claim. In particular, $\tau^rX=\tau^lY$, for some $l\in \Z$ and $Y\in \Sa$. Setting $s=r-l$, we see that the predecessors of $\tau^sX$ in $\GaA$ are all left stable. The proof of the lemma is completed.

\medskip

The following statement and its dual exhibit some interesting properties of not semi-stable modules in $\GaA$.

\medskip

\begin{Lemma}\label{nsp}

Let $M$ be a module in $\GaA$. If $M$ has infinitely many not left stable predecessors in $\GaA$, then all of its successors belong to an infinite right stable component of $\GaA$ which is not left stable and contains oriented cycles.

\end{Lemma}

\noindent{\it Proof.} Let $M_i$, $i\in \N$, be pairwise distinct and not left stable predecessors of $M$ in $\Ga$. Since $\GaA$ contains only finitely many projective modules, we may assume that there exists a projective module $P$ in $\GaA$ such that $M_i=\tau^{-n_i}P$ for some integer $n_i\ge 0$. Then, the $n_i$ with $i\ge 1$ are pairwise distinct, and in particular, they are unbounded. As a consequence, $P$ is right stable. By Lemma \ref{ss-module}(2), we obtain an integer $t>0$ such that the successors in $\GaA$ of $\tau^{-t}P$ belong to a right stable component $\Ga$, which is clearly not left stable.
Since the $n_i$ are unbounded, there exists no loss of generality in assuming that $n_i\ge t$ for all $i\ge 1$. In particular, $M$ is a successor of $\tau^{-t}P$. Therefore, the successors of $M$ in $\GaA$ are successors of $\tau^{-t}P$, and hence, they all belong to $\Ga$. Since each of the $\tau^{-n_i}P$ with $i\ge 1$ lies on a path from $\tau^{-t}P$ to $M$, by Lemma \ref{str-int-fin}, $\Ga$ contains oriented cycles. The proof of the lemma is completed.

\medskip

We are ready to given a number of combinatorial characterizations of almost acyclic components of $\GaA$. For this purpose, given a connected component $\mathcal{C}$ of $\GaA$, we define its {\it core} to be the full subquiver of $\cC$ generated by the modules lying on some path from a projective module to an injective module. Clearly, the core of $\mathcal{C}$ is convex in $\cC$.

\medskip

\begin{Theo}\label{thm-ifc}

Let $A$ be an artin algebra, and let $\mathcal{C}$ be a connected component of $\GaA$. The following conditions are equivalent.

\vspace{-2pt}

\begin{enumerate}[$(1)$]

\item The component $\mathcal{C}$ is almost acyclic.

\item The component $\mathcal{C}$ is interval-finite.

\item The core of $\hspace{0.4pt}\mathcal{C}$ is finite and contains all possible oriented cycles in $\mathcal{C}.$

\item Every infinite semi-stable component of $\mathcal{C}$ has no oriented cycle.

\end{enumerate}

\end{Theo}

\noindent{\it Proof.} First of all, it is evident that Statement (3) implies Statement (1). Suppose that $\mathcal{C}$ has an infinite semi-stable component $\Ga$ with oriented cycles. If $\Ga$ contains some $\tau$-periodic modules, then it is $\tau$-periodic. Being infinite, $\Ga$ is a stable tube; see, for example, \cite[(3.4)]{Liu0}. In this case, it is easy to see that $\Ga$ is not interval-finite and every module in $\Ga$ lies on an oriented cycle in $\Ga$. Assume that $\Ga$ contains no $\tau$-periodic module. In view of the second infinite path stated in Lemma \ref{ss-cpt}, we obtain a module $M\in \Ga$ and a positive integer $t$ such that $\Ga$ contains infinitely many oriented cycles $M\rightsquigarrow \tau ^{tn}M \rightsquigarrow M,$ where $n\in \Z$. In particular, $\Ga$ is not interval-finite and has infinitely many modules lying on oriented cycles. This shows that each of Statements (1) and (2) implies Statement (4).

It remains to show that Statement (4) implies Statements (2) and (3) hold. Suppose that Statement (4) holds. Assume on the contrary that $\mathcal{C}$ contains an infinite interval $[M, N]$ with $M, N$ some modules. That is, the full subquiver $\mathcal{D}$ of $\mathcal{C}$ generated by the modules lying on paths $M\rightsquigarrow N$ is infinite. By K\"{o}nig's Lemma, $\mathcal{D}$ has an infinite path
$$\xymatrix{(*) \quad \cdots \ar[r]  & N_i \ar[r] & \cdots \ar[r] & N_1 \ar[r] & N_0=N,}$$ where the $N_i$ with $i\ge 0$ are pairwise distinct. In view of Statement (4), we deduce from Lemma \ref{nsp} that there exists an integer $r$ such that $N_i$ is left stable for every $i\ge r$. Let $\Ga$ be the left stable component of $\mathcal{C}$ containing the $N_i$ with $i\ge r$. By Statement (4), $\Ga$ contains no oriented cycle, and by Lemma \ref{str-int-fin}, $\Ga$ is interval-finite. In particular, $M\not\in \Ga$. Then, for each $i\ge r$, there exists a path $\zeta_i: M\rightsquigarrow N_i$ in $\mathcal{D}$, which is the composite of a path $\,\xi_i: M\rightsquigarrow X_i,$ an arrow $\alpha_i: X_i\to Y_i$ with $X_i$ not left stable, and a path $\eta_i: Y_i\rightsquigarrow N_i$ in $\Ga$. By Lemma \ref{nsp}, the set $\{X_i \mid i\ge r\}$ is of finite cardinality, and so is $\{Y_i \mid i\ge r\}$. Therefore, we may assume that $Y_i=Y$ for some  $Y\in \Ga$ and all $i\ge r$. This yields infinitely many paths $\hspace{-1pt}\xymatrix{Y\ar@{~>}[r]^{\eta_i} & N_i\ar@{~>}[r]^{\omega_i} & N_r}\hspace{-2pt}$ in $\Ga$, where $\omega_i: N_i\rightsquigarrow N_r$ is the subpath of the infinite path $(*)$, a contradiction to Lemma \ref{str-int-fin}. This establishes Statement (2).
As a consequence, the core of $\mathcal{C}$, written as $\Oa$, is finite. Consider an oriented cycle
$$\sigma: \; \xymatrixcolsep{20pt}\xymatrix{Z_0\ar[r] & Z_1 \ar[r] & \cdots \ar[r] & Z_{n-1}\ar[r] & Z_n=Z_0}$$
in $\mathcal{C}$. Assume first that $\sigma$ is contained in a semi-stable component $\mathit\Theta$ of $\mathcal{C}$.
By Statement (4), $\mathit\Theta$ is finite, and hence, $\tau$-periodic.
It is well known that $\mathit\Theta \ne \mathcal{C}$; see \cite[(VII.2.1)]{ARS}. Thus, $\mathcal{C}$ contains an edge $U$---$\hspace{.4pt}V$ with $U\in \mathit\Theta$ and $V\not\in \mathit\Theta$. Since $U$ is $\tau$-periodic and $V$ is not, the $\tau$-orbit of $V$ contains a projective module $P$ and an injective module $I$. As a consequence, $\mathcal{C}$ has a path $P\rightsquigarrow U\rightsquigarrow I$. Let $Z\in \mathit\Theta$. Being $\tau$-periodic, $\mathit\Theta$ contains an oriented cycle $U\rightsquigarrow Z\rightsquigarrow U$, and hence, $Z$ lies in the core $\Oa$. In particular, $\sigma$ lies entirely in $\mathit\Omega$.

Assume next that $Z_s$ is not left stable and $Z_t$ is not right stable for some integers $0\le s, t\le n$. Then, $\mathcal{C}$ contains a path $P\rightsquigarrow Z_s$ with $P$ projective and a path $Z_t\rightsquigarrow I$ with $I$ injective. Since $\sigma$ is an oriented cycle, $\mathcal{C}$ contains paths $P\rightsquigarrow Z_i\rightsquigarrow I$
for all $1\le i\le n$. That is, $\sigma$ lies in $\mathit\Omega$. This proves Statement (3). The proof of the theorem is completed.

\medskip

\noindent{\sc Remark.} The almost acyclic components have been characterized by the existence of a multisection; see \cite[(2.5)]{RS}. Note that the core of a multisection of an almost acyclic component seems different from the core of the component defined here.

\medskip

\section{Components with bounded short cycles}

\medskip

The objective of this section is to study the connected components of $\GaA$ with bounded short cycles. This will yields a new characterization of representation-finite algebras, which includes a well known result of Ringle's saying that a representation-directed algebra is representation-finite; see \cite[(2.4)]{Rin}.

\medskip

\begin{Defn}

Let  $\Ga$ be a full subquiver of $\GaA$. A {\it cycle} in ${\rm rad}(\Ga)$ is a cycle in ${\rm ind}\hspace{.5pt}A$ passing through only modules in $\Ga$. We shall say that $\Ga$ is a subquiver with {\it bounded short cycles} if there exists a bound for the depths of all possibles short cycles in ${\rm rad}(\Ga)$, and otherwise, $\Ga$ is a subquiver with {\it unbounded short cycles}.

\end{Defn}

\medskip

Given any finite subquiver $\Ga$ of $\GaA$, it is evident that ${\rm rad}(\Ga)$ has bounded short cycle. The following result says that an infinite semi-stable component with bounded short cycles contains no oriented cycle.

\medskip

\begin{Lemma}\label{unbounded}

Let $\Ga$ be a semi-stable component of $\GaA$ with bounded short cycles. If $\Ga$ is infinite, then it is acyclic.

\end{Lemma}

\noindent{\it Proof.} We shall consider only the case where $\Ga$ is an infinite left stable component of $\GaA$. Suppose that $\Ga$ contains some $\tau$-periodic modules. Being infinite, $\Ga$ is a stable tube, say of rank $r$; see \cite{HPR}. Fix a module $X\in \Ga$. It is evident that there exists an infinite sectional path
$$\xymatrixcolsep{20pt}\xymatrix{X=X_0\ar[r] & X_1 \ar[r] & \cdots \ar[r] & X_n\ar[r] & \cdots}$$
in $\Ga$. Setting $Y_n=\tau^nX_n$ for all $n\ge 0$, we obtain another infinite sectional path
$$\xymatrix{\cdots \ar[r] & Y_n\ar[r] & \cdots \ar[r]& Y_1\ar[r] & Y_0=X} \vspace{1pt}$$ in $\Ga$. Observe that $Y_{rn}=\tau^{nr}X_{rn}=X_{rn}$ for all $n\ge 0$. Choosing irreducible maps $f_n: X_{n-1}\to X_n$ and $g_n: Y_n\to Y_{n-1}$ in ${\rm mod}A$ for all $n\ge 1$, we see that the maps $f_{nr}\cdots f_1: X_0\to X_{nr}$ and $g_{nr}\cdots g_1: X_{nr}\to X_0$ form a a short cycle of depth $nr\,;$ see \cite[(13.3)]{IgT}, a contradiction. Suppose now that $\Ga$ contains oriented cycles but no $\tau$-periodic module. By Proposition \ref{ss-cpt}, $\Ga$ contains two infinite sectional paths \vspace{-1pt}
$$\xymatrixcolsep{20pt}\xymatrix{
\cdots \ar[r] & \tau^{2t}M_1\ar[r] &  \tau^t M_r\ar[r] & \cdots \ar[r] & \tau^tM_1\ar[r] &
M_r \ar[r] & \cdots \ar[r] & M_1
}\vspace{-4pt}$$
and \vspace{-1pt} $$\xymatrixcolsep{20pt}\xymatrix{
N_1\ar[r]& \cdots \ar[r] & N_s\ar[r] & \tau^tN_1\ar[r] & \cdots \ar[r] & \tau^tN_s\ar[r] & \tau^{2t}N_1\ar[r] & \cdots
}$$
with $r, s, t>0$ and $M_1=N_1$. Given an integer $n$, using a similar argument as above, we can find a map $u_n: \tau^{\hspace{0.4pt}nt} M_1\to M_1$ of depth $rn$ and a map $v_n: M_1\to \tau^{\hspace{0.4pt}nt} M_1$ of depth $sn$. In particular, ${\rm rad}(\Ga)$ has short cycles of arbitrarily large depth, a contradiction. The proof of the lemma is completed.

\medskip

We shall regard a full subquiver of $\GaA$ as a translation quiver with the induced translation. The following result and its dual generalize slightly the result stated in  \cite[(2.1)]{Liu3}.

\medskip

\begin{Lemma}\label{predec-compl-sec}

Let $\Ga$ be a connected full subquiver of $\GaA$, and let $\Da$ be a section of $\Ga$ containing no right infinite path. If $\Ga$ contains all
the predecessors of $\Da$ in $\GaA$, then ${\rm ann}(\Da)={\rm ann}(\Oa)$, where $\Oa$ is the full subquiver of $\GaA$ generated by the predecessors of $\Da$ in $\GaA$.

\end{Lemma}

\noindent{\it Proof.} Containing no right infinite path, by Kr\"{o}nig's Lemma, $\Da$ contains only finite many paths starting at any given module. For each module $M\in \Da$, we shall denote by $d(M)$ the maximal length of the paths in $\Da$ starting at $M$.

Suppose that $\Oa$ is contained in $\Ga$. Write $I={\rm ann}(\Da)$. Fix a module $Y\in \Oa$. Then, $Y\in \tau^n\Da$ for some integer $n\ge 0$.
We shall show by an induction on $n$ that $IY=0$. Indeed, this is trivial for $n=0$. Suppose that $n>0$ and the statement holds for $n-1$. Then, $Y=\tau^nX$ with $X\in \Da$, and there exists in ${\rm mod}\hspace{0.5pt}A$ an almost split sequence
$$\xymatrixcolsep{20pt}\xymatrix{0\ar[r]& Y\ar[r]& Y_1\oplus \cdots \oplus Y_r \ar[r] & \tau^{n-1}X\ar[r] & 0,}$$
where $Y_1, \ldots, Y_r\in \Oa$. We claim that $IY_i=0$ for $i=1, \ldots, r$. Indeed,
for each $1\le i\le r$, there exists some $X_i\in \Da$ such that $Y_i=\tau^{n-1}X_i$ or $Y_i=\tau^nX_i$, where the second case occurs if and only if $\Da$ contains an arrow $X\to X_i$. We shall establish the claim by an induction on $d(X)$.
If $d(X)=0$, then $X$ is a sink vertex in $\Da$. Thus, $Y_i=\tau^{n-1}X_i$, and by the induction hypothesis on $n-1$, we obtain $IY_i=0$, for $i=1, \ldots, r$. Suppose now that $d(X)>0$. Fix an integer $i$ with $1\le i\le r$. If $Y_i=\tau^{n-1}X_i$, then $IY_i=0$ by the induction hypothesis on $n-1$. If $Y_i=\tau^nX_i$, then $\Da$ has an arrow $X\to X_i$, and thus,
$d(X_i)<d(X)$. Therefore, $IY_i=0$ by the induction hypothesis on $d(X)$. This establishes our claim. As a consequence, $IY=0$. The proof of the lemma is completed.

\medskip

We are now ready to obtain the first main result of this section.

\medskip

\begin{Theo}\label{scbc-ta}

Let $A$ be an artin algebra, and let $\Ga$ be an infinite left $($respectively, right$)$ stable component of $\GaA$. If ${\rm rad}(\Ga)$ has bounded short cycles, then $\Ga$ contains a finite section $\Da$ 
such that

\vspace{-2pt}

\begin{enumerate}[$(1)$]

\item $B=A/{\rm ann}(\Da)$ is a tilted algebra with $\Da$ being a slice of $\Ga_B;$

\vspace{1pt}

\item the predecessors $($respectively, successors$)$ of $\Da$ in $\GaA$ gene\-rate a predecessor-closed $($respectively, successor-closed $)$ subquiver of $\Ga_B$.

\end{enumerate}

\end{Theo}

\noindent{\it Proof.} Assume that ${\rm rad}(\Ga)$ has bounded short cycles. By Lemma \ref{unbounded}, $\Ga$ is acyclic, and hence, it contains a section $\Sa$ with a unique sink and no projective predecessor in $\GaA$; see \cite[(3.3)]{Liu1}. Now, $\Da=\tau \Sa$ is a section of $\Ga$ with a unique sink and no projective immediate successor in $\GaA$. Moreover, for any predecessor $X$ of $\Da$ in $\GaA$, both $X$ and $\tau^-X$ belong to $\Ga$. In particular, the full subquiver $\Oa$ of $\GaA$ generated by the predecessors of $\Da$ in $\GaA$ is contained in $\Ga$.

\vspace{1pt}

Let $\mathcal{C}$ be the connected component of $\GaA$ containing $\Ga$. Since $\Da$ has no immediate projective successor in $\GaA$, we see that $\Da$ is a cut of $\mathcal{C}$. If $\Da$ is infinite, then ${\rm rad}(\Ga)$ contains a short cycle
passing through predecessors of $\Da$ in $\Ga$; see the proofs of \cite[(2.2), (2.3), (2.4)]{Liu4}. Since $\Ga$ is acyclic and contains all the predecessors of $\Da$ in $\GaA$, this short cycle is of infinite depth, a contradiction. Thus, $\Da$ is finite.

We claim that $\rad^\infty(M, N)=0$ for all modules $M, N\in \Da$. Indeed, suppose on the contrary that $f_0: M\to N$ is a non-zero map in $\rad^\infty(M, N)$ with $M, N\in \Da$. Since $\Ga$ contains the predecessors of $\Da$ in $\GaA$, we can find an infinite path
$$\xymatrixcolsep{20pt}\xymatrix{\cdots\ar[r]& N_i\ar[r]& N_{i-1}\ar[r]& \cdots \ar[r] & N_1\ar[r]&N_0=N}$$ in $\Ga$
such that $\rad^\infty(M, N_i)$ has a non-zero map $f_i$, for every $i\ge 0$; see \cite[(1.1)]{Liu5}. Since $\Da$ is a finite section of $\Ga$, there exists a minimal $s\ge 0$ such that $\Ga$ has a path
$$\xymatrixcolsep{20pt}\xymatrix{\rho: & N_s=M_t\ar[r]& M_{t-1}\ar[r]& \cdots \ar[r] & M_1\ar[r]&M_0=M.}$$

Suppose that $\rho$ is not sectional with $M_{j+1}=\tau M_{j-1}$ for some $0< j<t-1$. Since $\Da$ is a section of $\Ga$, we have $s>0$. By our choice of $\Da$, the $\tau^-M_i$ with $t\ge i \ge j+1$ belong to $\Ga$. This yields a path
$$\xymatrixcolsep{20pt}\xymatrix{N_{s-1}\ar[r]& \taum M_t\ar[r]& \cdots \ar[r] & \taum M_{j+2}\ar[r] & M_{j-1} \ar[r]& \cdots \ar[r] &M_0=M,}$$
a contradiction to the minimality of $s$. Therefore, $\rho$ is sectional, and consequently, $\Hom_A(N_s, M)\ne 0$. Thus, we obtain a short cycle $\hspace{-2pt}\xymatrixcolsep{16pt}\xymatrix{M\ar[r]^{f_s} & N_s\ar[r]^g & M \hspace{-2pt}}$ of infinite depth in ${\rm rad}(\Ga)$, a contradiction. This establishes our claim.

Next, suppose that there exists a non-zero map $f_0: M\to \tau N_0$ with $M, N_0\in \Da.$ Since $\tau N_0$ is not a predecessor of $M$ in $\Ga$, it is not a predecessor of $M$ in $\GaA$. Hence, $f_0\in \rad^\infty(M, \tau N_0)$. Considering the minimal left almost split monomorphism for $\tau N_0$, we obtain an irreducible map $f_1: \tau N_0\to M_1$ with $M_1\in \Ga$ such that $0\ne f_1f_0\in \rad^\infty(M, M_1).$ By the above claim, $M_1\not\in \Da$, and thus, $M_1=\tau N_1$ with $N_1\in \Da$. This yields an arrow $N_0\to N_1$ in $\Da$ with $\rad^\infty(M, \tau N_1)\ne 0$. Continuing this process, we obtain an infinite path
$$\xymatrixcolsep{20pt}\xymatrix{N_0\ar[r]& N_1 \ar[r] & \cdots \ar[r] & N_{r-1} \ar[r]& N_r\ar[r]& \cdots}$$
in $\Da$, contrary to the finiteness of $\Da$. This shows that $\Hom_A(\Da, \tau \Da)=0$. As a consequence, $B=A/{\rm ann}(\Da)$ is a tilted algebra with $\Da$ being a slice of $\Ga_B$; see \cite[(2.8)]{Liu5}. By Lemma \ref{predec-compl-sec}, the predecessors of $\Da$ in $\GaA$ are $B$-modules. Therefore, $\Omega$ is a subquiver of $\Ga_B$. Since $\Oa$ is left stable and predecessor-closed in $\GaA$, we see that $\Oa$ is predecessor-closed in $\Ga_B.$ The proof of the proposition is completed.

\medskip

The following consequence of Theorem \ref{scbc-ta} includes the result on semi-regular components admitting no short cycles obtained in \cite[(2.7)]{Liu4}.

\medskip

\begin{Theo}\label{semi-reg}

Let $A$ be an artin algebra. If $\mathcal{C}$ is a semi-regular component of $\GaA$, then ${\rm rad}(\mathcal{C})$ has bounded short cycles if and only if $B=A/{\rm ann}(\mathcal{C})$ is tilted with $\mathcal C$ a connecting component of $\Ga_B$.

\end{Theo}

\noindent{\it Proof.} The sufficiency is evident, see, for example, \cite[(2.7)]{Liu4}. Let $\mathcal{C}$ be a connected component of $\GaA$ with no projective module. In particular, $\mathcal{C}$ is infinite. Assume that ${\rm rad}(\Ga)$ has bounded short cycles. By Theorem \ref{scbc-ta}, there exists a section $\Da$ of $\mathcal{C}$ such that $B=A/{\rm ann}(\Da)$ is a tilted algebra with $\mathcal C$ a connecting component of $\Ga_B$. Since ${\rm ann}(\Da)={\rm ann}(\mathcal{C})$; see \cite[(2.1)]{Liu4}, we obtain $B=A/{\rm ann}(\mathcal{C})$. The proof of the theorem is complete.

\medskip

The following statement extends some results stated in \cite[(2.6),(2.8)]{Liu4} on Auslander-Reiten components admitting short cycles.

\medskip

\begin{Theo}\label{almost_acyclic}

Let $A$ be an artin algebra. The Auslander-Reiten quiver $\GaA$ has at most finitely many connected components with bounded short cycles, and each of them is almost acyclic with only finitely many $\tau$-orbits.

\end{Theo}

\noindent{\it Proof.} Let $\mathcal{C}$ be a connected component of $\GaA$ with bounded short cycles. Observe that $\mathcal{C}$ has only finitely many non-trivial semi-stable components, and their union is co-finite in $\cC$; see \cite[(3.1)]{Liu1}. Having a finite section by Theorem \ref{scbc-ta}, every infinite semi-stable component of $\mathcal{C}$ is acyclic with only finitely many $\tau$-orbits. Thus, $\mathcal{C}$ has only finitely many $\tau$-orbits, and by Theorem \ref{thm-ifc}, it is almost acyclic.

If $\mathcal{C}$ is semi-regular then, by Proposition \ref{semi-reg}, $\mathcal C$ is a connecting component of a tilted algebra $A/{\rm ann}(\mathcal{C})$. In particular, ${\rm rad}(\mathcal{C})$ contains no short cycle; see \cite[(2.7)]{Liu4}. Having at most finitely many connected components which are not semi-regular and at most finitely many connected components with no short cycle; see \cite[(2.8)]{Liu4}, $\GaA$ has at most finitely many connected components with bounded short cycles. The proof of the theorem is completed.

\medskip

A connected artin algebra $A$ is {\it generalized double tilted} if and only if $\GaA$ contains a faithful, almost acyclic and generalized standard component, which is called a {\it connecting component}; see \cite[Section 3]{RS}.

\medskip

\begin{Prop}\label{main-th-3}

Let $A$ be an artin algebra. A connected component $\mathcal{C}$ of $\GaA$ is generalized standard with bounded short cycles if and only if $B=A/{\rm ann}(\mathcal C)$ is generalized double tilted with $\mathcal C$ a connecting component of $\Ga_B$.

\end{Prop}

\noindent{\it Proof.} The necessity follows immediately from Theorem \ref{almost_acyclic}. Assume now that $B=A/{\rm ann}(\mathcal C)$ is generalized double tilted and $\mathcal C$ is a connecting component of $\it\Gamma_B$. Being almost acyclic, by Theorem \ref{thm-ifc},
$\mathcal{C}$ contains a finite core $\it\Omega$. Let $b$ be the maximal $R$-length of the modules in $\it\Omega$.
Consider a short cycle $\sigma$ consisting of two maps $f: M\to N$ and $g: N\to M$ in ${\rm rad}(\mathcal{C})$. Since $\mathcal{C}$ is generalized standard, we deduce from Theorem \ref{thm-ifc} that both $f$ and $g$ are sums of composites of irreducible maps in ${\rm rad}(\Oa)$. In view of the Harada-Sai Lemma; see \cite{HS}, and also \cite[(VI.1.3)]{ARS}, we see that both $f$ and $g$ are of depth less than $2^b$. That is, ${\rm dp}(\sigma)<2^b.$ The proof of the proposition is completed.

\medskip

\noindent{\sc Example.} Let $A=kQ/I$, where $k$ is a field, $Q$ is the quiver

\vspace{-5pt}

$$\xymatrix{5\ar[r]&4\ar@<0.5ex>[r]\ar@<-0.5ex>[r] & 3\ar[r] & 2\ar@<0.5ex>[r] &1,\ar@<0.5ex>[l]}$$
and $I$ is the ideal in the path algebra $kQ$ generated by the paths of length two. It is easy to see that $\Ga_A$ contains a generalized standard component with bounded short cycles as follows$\,:$
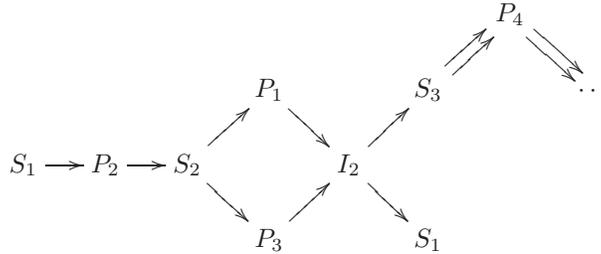
$$\xymatrix@R=0.5cm@C=0.5cm{&&&&&&P_4\ar@<0.5ex>[rd]\ar@<-0.5ex>[rd] \\&&&P_1\ar[rd]&&S_3\ar@<0.5ex>[ru]\ar@<-0.5ex>[ru]&&\cdots \\S_1\ar[r]&P_2\ar[r]&S_2\ar[ru]\ar[rd]&&I_2\ar[rd]\ar[ru]\\
&&&P_3\ar[ru]&&S_1}$$

\medskip

\medskip

It has been shown that the artin algebra $A$ is representation-finite if ${\rm mod}\hspace{0.5pt}A$ contains no short cycle; see \cite{HL}.
In order to improve this result, we shall say that ${\rm mod}\hspace{.5pt}A$ has {\it bounded short cycles} if there exists a bound for the depths of the maps on short cycles in $\ind\hspace{.5pt}A$.

\medskip

\begin{Theo}\label{last-thm}

An artin algebra $A$ is of finite representation type if and only if ${\rm mod}\hspace{.5pt}A$ has bounded short cycles.

\end{Theo}

\noindent{\it Proof.} 
Suppose first that $A$ is of finite representation type. Then ${\rm rad}^n({\rm mod}\hspace{.5pt}A)=0$ for some integer $n>0$; see \cite[(V.7.7)]{ARS}. That is, every non-zero map in ${\rm ind}\hspace{.5pt}A$ is of depth less than $n$. In particular, ${\rm mod}\hspace{.5pt}A$ has bounded short cycles.

Suppose now that ${\rm mod}\hspace{.5pt}A$ has bounded short cycles but is of infinite type. Then, $\GaA$ has an infinite connected component
$\mathcal C$; see \cite[(VI.1.4)]{ARS}. We may assume with no loss of generality that $\mathcal{C}$ contains an infinite left-stable component $\Ga$; see \cite[(3.1)]{Liu1}. By Theorem \ref{scbc-ta}, $\Ga$ contains a section $\Da$ such that $B=A/{\rm ann}(\Da)$ is a tilted algebra, and the predecessors of $\Da$ in $\GaA$ generate a predecessor-closed subquiver of a connecting component of $\Ga_B$. In particular, $B$ is of representation-infinite. Then, $\Ga_B$ has a connected component $\mathscr{C}$ containing non-directing modules; see \cite[(2.4)]{Rin1}, and also \cite{HL}. Observe that $\mathscr{C}$ cannot be a prepropjective component, a preinjective component or a connecting component of $\Ga_B$; see \cite[(2.4), (4.2)]{Rin1}. Therefore, $\mathscr{C}$ is either quasi-serial or is obtained from a quasi-serial translation quiver by ray insertions or by co-ray insertions; see  \cite[(3.7)]{Liu6}. If $\mathscr{C}$ contains oriented cycles, by Lemma \ref{unbounded}, ${\rm rad}(\mathscr{C})$ has unbounded short cycles in ${\rm ind}\hspace{.5pt}B$. Otherwise, $\mathscr{C}$ is obtained from a translation quiver of shape $\Z \A_\infty$ by ray insertions or by co-ray insertions. In particular, $\mathscr{C}$ has infinitely many $\tau_B$-orbits, and by Theorem \ref{scbc-ta}, ${\rm rad}(\mathscr{C})$ has unbounded short cycles. In all cases, ${\rm mod}\hspace{.5pt}B$ has unbounded short cycles, and so does ${\rm mod}\hspace{.5pt}A$, a contradiction. The proof of the theorem is completed.

\bigskip


\begin{thebibliography}{99}

\bigskip

%

%

\bibitem{AR} {\sc M. Auslander, I. Reiten}, ``Modules determined by their composition factors," Illinois J. Math. 29 (1985) 280 - 301.

\smallskip

\bibitem{ARS} {\sc M. Auslander, I. Reiten and S. Smal\o}, ``Representation Theory of Artin Algebras," Cambridge
Studies in Advanced Mathematics 36 (Cambridge University Press, Cambridge, 1995).

%

%
%
%

%

%

%

\smallskip

\bibitem{HL} {\sc D. Happel and S. Liu}, ``Module categories without short cycles are of finite type," Proc. Amer. Math. Soc. 120 (1994) 371 - 375.

\smallskip

\bibitem{HPR}  {\sc D. Happel, U. Preiser and C. M. Ringel}, ``Vinberg's characterization of Dynkin diagrams using subadditive functions with application to DTr-periodic modules," Lecture Notes in Math. 832 (Springer, Berlin, 1980) 280 - 294.

\smallskip

\bibitem{HaR} {\sc D. Happel and C. M. Ringel,} ``Tilted algebras," Trans. Amer. Math. Soc. 274 (1982) 399 - 443.

\smallskip

\bibitem{HS} {\sc M. Harada and Y. Sai}, ``On categories of indecomposable modules {\rm I}," Osaka J. Math. 7 (1970) 323 - 344.

\smallskip

\bibitem{IgT} {\sc K. Igusa and G. Todorov,} ``A characterization of finite Auslander-Reiten quivers," J. Algebra 89 (1984) 148 - 177.


\smallskip

\bibitem{Liu} {\sc S. Liu}, ``Degrees of irreducible maps and the shapes of Auslander-Reiten quivers," J. London Math.
Soc. 45 (1992) 32 - 54.

\smallskip

\bibitem{Liu0} {\sc S. Liu}, ``Semi-stable components of an Auslander-Reiten quiver," J. London Math. Soc. 47 (1993) 405 - 416.

\smallskip

\bibitem{Liu3} {\sc S. Liu}, ``Tilted algebras and generalized standard Auslander-Reiten components," Arch. Math. 61 (1993) 12 - 19.

\smallskip

\bibitem{Liu6} {\sc S. Liu}, ``The connected components of the Auslander-Reiten quiver of a tilted algebra," J. Algebra 161 (1993) 505 - 523.

\smallskip

\bibitem{Liu1} {\sc S. Liu}, ``Shapes of connected components of the Auslander-Reiten quivers of artin algebras," Canad. Math. Soc. Conf. Proc. 19 (1995) 109 - 137.

\smallskip

\bibitem{Liu4} {\sc S. Liu}, ``On short cycles in a module category", J. London math. Soc. 51 (1995) 62 - 74.

\smallskip

\bibitem{Liu2} {\sc S. Liu}, ``Auslander-Reiten theory in a Krull-Schmidt category," Sao Paulo J. Math. Sci. 4 (2010) 425 - 472.

\smallskip

\bibitem{Liu5} {\sc S. Liu}, ``Another characterization of tilted algebras," Arch. Math. 104 (2015) 111 - 123.

%

\smallskip

\bibitem{MPS1} {\sc P. Malicki, J. A. de la Pe\~na and A. Skowro\'nski}, ``Cycle-finite module categories,"
{\it in}: Algebras, Quivers and Representations, Abel Symp. 8 (Springer, Verlag, 2013) 209 - 252.

\smallskip

\bibitem{MPS2} {\sc P. Malicki, J.A. de la Pe\~na and A. Skowro\'nski}, ``Finite cycles of indecomposable modules,"
 J. Pure Appl. Algebra 219 (2015) 1761 - 1799.

\smallskip

\bibitem{PeX} {\sc L. Peng and J. Xiao}, ``On the number of DTr-orbits containing directing modules", Proc. Amer. Math. Soc.
118 (1993) 753 - 756.

\smallskip

\bibitem{RS} {\sc I. Reiten and A. Skowro\'nski}, ``Generalized double tilted algebras," J. Math. Soc. Japan 56 (2004) 269 - 288.

\smallskip

\bibitem{RSS1} {\sc I. Reiten, A. Skowro\'nski and S. O. Smal\o}, ``Short chains and short cycles of modules," Proc. Amer. Math. Soc. 117 (1993) 343-354.

%
%
%

\smallskip

\bibitem{Rin} {\sc C. M. Ringel}, ``Tame Algebras and Integral Quadratic Forms," Lecture Notes in Mathematics 1099 (Springer, Verlag, Berlin, 1984).

\smallskip

\bibitem{Rin1} {\sc C. M. Ringel}, ``Representation theory of finite dimensional algebras," London Math. Soc. Lecture Note Ser. 116 (1986) 7 - 79.

\smallskip

\bibitem{Sky} {\sc A. Skowyrski},  ``A characterization of cycle-finite generalized double tilted algebras," J. Algebra 416 (2014) 1 - 24.

\smallskip

\bibitem{Sko} {\sc A. Skowro\'nski}, ``Cycles in module categories", NATO Adv. Sci. Inst. Ser. C Math. Phys. Sci. 424 (Kluwer Acad. Publ., Dordrecht, 1994) 309 - 345.

%

%

\smallskip

\bibitem{Sko3} {\sc A. Skowro\'nski}, ``Generalized standard Auslander-Reiten components," J. Math. Soc. Japan 46 (1994) 517 - 543.

\smallskip

\bibitem{Sko4} {\sc A. Skowro\'nski}, ``Cycle-finite algebras," J.  Pure Appl. Algebra 103 (1995) 105 - 116.

\smallskip

\bibitem{SkS} {\sc A. Skowro\'nski and S. O.  Smal\o}, ``Directing modules", J. Algebra 147 (1992) 137 - 146.

%

\bigskip

\end{thebibliography}
\end{document}